\documentclass[10pt]{article}
\usepackage{amssymb,amsmath,textcomp}
\newtheorem{satz}{Theorem}[section]
\newtheorem{lemma}{Lemma}[section]

\newtheorem{prop}{Proposition}[section]
\newcommand{\iR}{\mathbb{R}}
\newcommand{\iN}{\mathbb{N}}

\newcommand{\oH}{\hspace*{0.39em}\raisebox{0.6ex}{\textdegree}\hspace{-0.72em}H}
\DeclareMathOperator*{\esup}{ess\,sup}
\DeclareMathOperator*{\einf}{ess\,inf}

\begin{document}
\begin{center}
{\bf\Large A priori bounds for degenerate and singular evolutionary
partial integro-differential equations}
\end{center}
\vspace{0.7em}
\begin{center}
Vicente Vergara and Rico Zacher
\end{center}
\footnotetext{The first author was partially supported by FONDECYT,
project number 11070003, and by Convenio de Desempe\~no Universidad
de Tarapac\'a - Mineduc. The second author thanks FONDECYT for the
financial support given when he visited the first author in Arica
(Chile) in 2009.} \vspace{0.7em}
\begin{abstract}
We study quasilinear evolutionary partial integro-differential
equations of second order which include time fractional $p$-Laplace
equations of time order less than one. By means of suitable energy
estimates and De Giorgi's iteration technique we establish results
asserting the global boundedness of appropriately defined weak
solutions of these problems. We also show that a maximum principle
is valid for such equations.
\end{abstract}
\vspace{0.7em}
\begin{center}
\end{center}

\section{Introduction and main result}
Let $T>0$, and $\Omega$ be a bounded domain in $\iR^N$. In this
paper we are concerned with global a priori bounds for weak
solutions of quasilinear problems of the form
\begin{equation} \label{MProb}
\partial_t \Big(k\ast(u-u_0)\Big)-\mbox{div}\,a(t,x,u,Du)=b(t,x,u,Du),\;\;t\in (0,T),\,x\in
\Omega,
\end{equation}
where $Du$ stands for the gradient of $u$ w.r.t.\ the spatial
variables, $k\in L_{1,\,loc}(\iR_+)$ is a singular kernel, and
$k\ast v$ denotes the convolution on the positive halfline w.r.t.\
the time variable, that is $(k\ast v)(t)=\int_0^t
k(t-\tau)v(\tau)\,d\tau$, $t\ge 0$.

We will assume that the kernel $k$ satisfies the following
conditions.
\begin{itemize}
\item[{\bf (K1)}] $k$ is of type ${\cal PC}$, that is
(cf. \cite{ZWH},\cite{Za}) $k\in L_{1,\,loc}(\iR_+)$ is nonnegative
and nonincreasing, and there exists a kernel $l\in
L_{1,\,loc}(\iR_+)$ such that $k\ast l=1$ in $(0,\infty)$.
\item[{\bf (K2)}] $l\in L_q([0,T])$ for some $q>1$.
\end{itemize}

An important example is given by
\begin{equation} \label{paar}
k(t)=g_{1-\alpha}(t)e^{-\mu t}\quad\mbox{and}\quad
l(t)=g_{\alpha}(t)e^{-\mu
t}+\mu(1\ast[g_{\alpha}(\cdot)e^{-\mu\cdot}])(t), \quad t>0,
\end{equation}
with $\alpha\in (0,1)$ and $\mu\ge 0$, see also
\cite{ZWH},\cite{Za}. Here $g_\beta$ denotes the Riemann-Liouville
kernel
\begin{equation} \label{RLkernel}
g_\beta(t)=\,\frac{t^{\beta-1}}{\Gamma(\beta)}\,,\quad
t>0,\;\beta>0.
\end{equation}
In this case, (\ref{MProb}) amounts to a time fractional equation of
order $\alpha\in (0,1)$. Recall that for a (sufficiently smooth)
function $v$ on $\iR_+$, the Riemann-Liouville fractional derivative
$D_t^\alpha v$ of order $\alpha\in (0,1)$ is defined by $D_t^\alpha
v=\,\frac{d}{dt}\,(g_{1-\alpha}\ast v)$.

Letting $p>1$ and $\Omega_T=(0,T)\times \Omega$, we will further
assume that the functions $a:\Omega_T\times \iR^{N+1}\rightarrow
\iR^N$ and $b:\Omega_T\times \iR^{N+1}\rightarrow \iR$ are
measurable and that they satisfy the structure conditions
\begin{itemize}
\item[{\bf (Q1)}] $\quad\quad\quad\quad \quad\quad\quad(a(t,x,\xi,\eta)|\eta) \ge \,
C_0|\eta|^p-c_0|\xi|^\gamma-\varphi_0(t,x)$,
\item[{\bf (Q2)}] $\quad\quad\quad\quad\quad\quad\quad |a(t,x,\xi,\eta)| \,\le \,
C_1|\eta|^{p-1}+c_1|\xi|^{r\frac{p-1}{p}}+\varphi_1(t,x)$,
\item[{\bf (Q3)}] $\quad\quad\quad\quad\quad\quad\quad |b(t,x,\xi,\eta)| \,\le \,
C_2|\eta|^{p\frac{\gamma-1}{\gamma}}+c_2|\xi|^{\gamma-1}+\varphi_2(t,x)$,
\end{itemize}
for a.a. $(t,x)\in \Omega_T$, and all $\xi\in\iR$, $\eta\in \iR^N$.
Here $C_i, c_i$, $i=0,1,2$, are positive constants, and
\begin{itemize}
\item[{\bf (Q4)}] The parameter $\gamma$ lies in the range
\[
1<\gamma<
\frac{1-\frac{1}{q}+\frac{2}{N}}{\frac{1}{p}\big(1-\frac{1}{q}\big)+\frac{1}{Nq}}=:r.
\]
\item[{\bf (Q5)}] The functions $\varphi_i$, $i=0,1,2$, defined on
$\Omega_T$ are nonnegative, $\varphi_1^\frac{p}{p-1}\in
L_1(\Omega_T)$, and $\varphi_0, \,\varphi_2\in L_{s}(\Omega_T)$,
where
\begin{equation} \label{scond}
s>\frac{\frac{1}{p}\big(1-\frac{1}{q}\big)+\frac{1}{N}}{\frac{1}{N}\big(1-\frac{1}{q}\big)}=\,\frac{N}{p}\,+q'.
\end{equation}
\end{itemize}
Here, as usual, $q$ and $q'$ denote conjugate exponents, i.e.
$\frac{1}{q}+\frac{1}{q'}=1$.

The function $u_0=u_0(x)$ is a given data and plays the role of the
initial data for the function $u$. We will assume that $u_0\in
L_2(\Omega)$.

Throughout the paper we will further assume that $\partial \Omega$
satisfies the property of positive density, see Section \ref{SecP}.

Before describing the main results we give some comments on
applications. Problems of the form (\ref{MProb}) arise for example
in mathematical physics when describing dynamic processes in
materials with memory, e.g.\ in the theory of heat conduction with
memory, see \cite{JanI} and the references therein. Time fractional
diffusion equations which are obtained by taking $k=g_{1-\alpha}$ in
(\ref{MProb}) are also used to model anomalous diffusion, see e.g.\
\cite{Metz}. In this context, these equations are termed {\em
subdiffusion equations} (the time order $\alpha$ lies in $(0,1)$);
in the case $\alpha\in(1,2)$, which is not considered here, one
speaks of {\em superdiffusion equations}. We point out that our
general setting also includes models which describe nonlinear
diffusion phenomena. An important special case of (\ref{MProb}) is
the class of time fractional $p$-Laplace equations like e.g.\
(\ref{plaplace}) below. Let us further mention that time fractional
diffusion equations of time order $\alpha\in (0,1)$ are closely
related to a class of Montroll-Weiss continuous time random walk
models where the waiting time density behaves as $ t^{-\alpha-1}$
for $t\to \infty$, see e.g. \cite{Hilfer1}, \cite{Hilfer2},
\cite{Metz}.

We say that a function $u$ is a {\bf weak solution (subsolution,
supersolution)} of (\ref{MProb}) in $\Omega_T$, if $u$ belongs to
the space
\begin{align*}
\tilde{V}_{q,p}:=\{&\,v\in L_{2q}([0,T];L_2(\Omega))\cap
L_p([0,T];H^1_p(\Omega))\;
\mbox{such that}\;\\
&\;\;k\ast v\in C([0,T];L_2(\Omega)), \;\mbox{and}\;(k\ast
v)|_{t=0}=0\},
\end{align*}
$a(t,x,u,Du)$ and $b(t,x,u,Du)$ are measurable, and for any
nonnegative test function
\[
\eta\in \oH^{1,1}_2(\Omega_T):=H^1_2([0,T];L_2(\Omega))\cap
L_p([0,T];\oH^1_p(\Omega)) \quad\quad
\Big(\oH^1_p(\Omega):=\overline{C_0^\infty(\Omega)}\,{}^{H^1_p(\Omega)}\Big)
\]
with $\eta|_{t=T}=0$ there holds
\begin{equation} \label{WQ}
\int_{0}^{T} \int_\Omega \Big(-\eta_t [k\ast (u-u_0)]+
\big(a(t,x,u,Du)|D
\eta\big)-b(t,x,u,Du)\eta\Big)\,dx\,dt=(\le,\,\ge)\,\,0.
\end{equation}
This definition makes sense, since under conditions (Q1)-(Q5) the
integral in (\ref{WQ}) is finite, by H\"older's inequality and the
parabolic embedding $\tilde{V}_{q,p}\hookrightarrow L_r(\Omega_T)$
(see Proposition \ref{parabolicembed} below). We point out that
(\ref{MProb}) is considered without any boundary conditions, in this
sense weak solutions of (\ref{MProb}) as defined above are {\em
local} ones w.r.t.\ space. We further remark that weak solutions of
(\ref{MProb}) in the class $\tilde{V}_{q,p}$ have been constructed
in \cite{ZWH} in the linear case with $p=2$. In view of the basic
energy estimate (see below) and the known results in the case $p=2$
the space $\tilde{V}_{q,p}$ is the natural choice for weak solutions
in the general case $p\in (1,\infty)$. We strongly believe that
under stronger assumptions on the nonlinearities $a$ and $b$ it is
possible to prove the existence of weak solutions of (\ref{MProb})
in the class $\tilde{V}_{q,p}$ by means of the theory of monotone
operators and the techniques developed in \cite{ZWH}. Notice also
that the initial condition $u|_{t=0}=u_0$ has to be understood in a
weak sense. One can show (\cite{ZWH}) that in case of sufficiently
smooth functions $u$ and $k\ast(u-u_0)$, the condition $(k\ast
u)|_{t=0}=0$ implies $u|_{t=0}=u_0$.

To state our main results we set $\Gamma_T=(0,T)\times
\partial \Omega$ and $y_+:=\max\{y,0\}$. By $|A|$ we denote the Lebesgue measure of a measurable set
$A\subset \iR^N$. When we say that a function $u\in \tilde{V}_{q,p}$
satisfies $u\le K$ a.e.\ on $\Gamma_T$ for some number $K\in \iR$ we
mean that $(u-K)_+\in L_p([0,T];\oH^1_p(\Omega))$, likewise for
lower bounds on $\Gamma_T$. This convention allows to formulate our
results without extra smoothness assumptions on the boundary
$\partial\Omega$. Our main result reads as follows.
\begin{satz} \label{mainresult}
Let $p>1$, $T>0$, and $\Omega\subset \iR^N$ be a bounded domain. Let
the assumptions (K1),(K2),(Q1)-(Q5) be satisfied.

(i) (Subsolutions) Suppose $u_0\in L_2(\Omega)$ and that $K\ge 0$ is
such that $u_0\le K$ a.e.\ in $\Omega$. Then there exists a constant
\begin{equation} \label{constant}
C=C(N,p,q,C_0,c_0,C_2,c_2,\gamma,s,|l|_{L_q([0,T])},|\varphi_0+\varphi_2|_{L_s(\Omega_T)},T,|\Omega|)
\end{equation}
such that for any weak subsolution $u\in\tilde{V}_{q,p}$ of
(\ref{MProb}) in $\Omega_T$ satisfying $u\le K$ a.e.\ on $\Gamma_T$
there holds
\begin{equation} \label{mainbound0}
\esup_{\Omega_T} u\le 2\Big(K+\max\Big\{1,C\big(\int_0^T
\int_{\Omega}
u_+^\gamma\,dx\,dt\big)^{\frac{\theta}{r-\gamma}}\Big\}\Big),
\end{equation}
where $r$ is defined in (Q4) and
\[
\theta=\frac{\frac{1}{N}\big(1-\frac{1}{q}\big)-\frac{1}{s}\big[\frac{1}{p}\big(1-\frac{1}{q}\big)+\frac{1}{N}\big]}
{\frac{1}{p}\big(1-\frac{1}{q}\big)+\frac{1}{Nq}}.
\]
(ii) (Supersolutions) Suppose $u_0\in L_2(\Omega)$ and that $K\ge 0$
is such that $u_0\ge -K$ a.e.\ in $\Omega$. Then there exists a
constant $C$ like in (\ref{constant}) such that for any weak
supersolution $u\in\tilde{V}_{q,p}$ of (\ref{MProb}) in $\Omega_T$
satisfying $u\ge -K$ a.e.\ on $\Gamma_T$ there holds
\begin{equation} \label{mainbound2}
\einf_{\Omega_T} u\ge -2\Big(K+\max\Big\{1,C\big(\int_0^T
\int_{\Omega}
(-u)_+^\gamma\,dx\,dt\big)^{\frac{\theta}{r-\gamma}}\Big\}\Big).
\end{equation}
\end{satz}
Note that $c_1$, $C_1$, and $\varphi_1$, which appear in (Q2), do
not play any role in determining the constant in (\ref{mainbound0})
and (\ref{mainbound2}), respectively.

An important special case of (\ref{MProb}) is the equation
\begin{equation} \label{plaplace}
\partial_t^\alpha (u-u_0)-\mbox{div}\,\big(|Du|^{p-2}Du\big)=f\quad
\mbox{in}\;\Omega_T,
\end{equation}
with $\alpha\in (0,1)$ and $p>1$. We have the following result.
\begin{satz} \label{plaplacebound}
Let $\alpha\in (0,1)$, $p>1$, $T>0$, and $\Omega\subset \iR^N$ be a
bounded domain. Suppose that $u_0\in L_\infty(\Omega)$ and that
$f\in L_s(\Omega_T)$ with $s>\frac{N}{p}+\frac{1}{\alpha}$. Let
further $q>1$ be a fixed number satisfying
\begin{equation} \label{qcond}
s>\frac{N}{p}+q'>\frac{N}{p}+\frac{1}{\alpha}.
\end{equation}
Then for any weak solution $u\in\tilde{V}_{q,p}$ of (\ref{plaplace})
in $\Omega_T$ which is essentially bounded on $\Gamma_T$ there holds
\begin{equation} \label{alphapbound}
|u|_{L_\infty(\Omega_T)}\le
C\big(N,p,\alpha,s,|f|_{L_s(\Omega_T)},T,|\Omega|,\max\{|u_0|_{L_\infty(\Omega)},\esup_{\Gamma_T}
|u|\}\big).
\end{equation}
\end{satz}
Here the condition on $f$ is sharp, at least in the cases $p=2$,
$\alpha\in(0,1)$ and $p>2$, $\alpha\in(p'/N,1)$, as we will show in
Section \ref{pl}. We also remark that the estimate
(\ref{alphapbound}) is stable for $\alpha\rightarrow 1$, that is,
the constants in the proof remain bounded as $\alpha\rightarrow 1$.
Hence in this sense we recover well-known results for equations like
the classical parabolic $p$-Laplace equation, which can be found in
the monograph \cite{DB}.

In this paper we further prove that in case of so-called homogenous
structures (see Section \ref{SecHomStruc}) the weak maximum
principle for weak solutions takes the same form as in the classical
parabolic case. This applies e.g.\ to equation (\ref{plaplace}) with
$f=0$.

In the literature not much seems to be known concerning a regularity
theory for {\em weak} solutions to (\ref{MProb}) in the general
setting considered in this paper. To our knowledge, the only paper
in this direction is \cite{Za}, where the global boundedness of weak
solutions was proved in the case $p=2$ under similar assumptions on
the kernel $k$ and the nonlinearities $a$ and $b$. On the other hand
there exists a rather well developed regularity theory for
degenerate ($p>2$) and singular ($1<p<2$) parabolic equations of the
form (\ref{MProb}) with $\partial_t (k\ast(u-u_0))$ replaced by
$\partial_t u$, see the monograph \cite{DB} and the references given
therein as well as the recent work \cite{DBGV}. This theory includes
besides local and global $L_\infty$-bounds also much deeper results
such as Harnack and H\"older estimates for weak solutions. For the
case $p=2$ we also refer to \cite{LSU} and \cite{Lm}. In the time
fractional case the situation is much harder due to the nonlocal
nature of $\partial_t^\alpha$. Recently, a weak Harnack inequality
was proved for nonnegative weak supersolutions of (\ref{MProb}) in a
special case where $p=2$ and $k=g_{1-\alpha}$, see \cite{Za2}.
Concerning results in {\em stronger} settings for (\ref{MProb}) as
well as abstract variants of it (mostly with $p=2$) we refer to
\cite{Ba}, \cite{CLS}, \cite{CP2}, \cite{Koch}, \cite{Grip1},
\cite{JanI}, \cite{ZEQ}, \cite{ZQ}.

Our proofs of the global $L_\infty$-bounds use De Giorgi's iteration
technique and are based on suitable truncated energy estimates for
weak solutions of (\ref{MProb}). These estimates are derived by
combining the techniques from \cite{Za} and \cite{DB}. A key
ingredient is the basic inequality (\ref{BI}) (see below) for
nonnegative nonincreasing kernels. We further adopt the method of
time regularization of the equation which goes back to \cite{Za} in
the weak setting (see also \cite{VZ}) and uses the Yosida
approximations of the operator $B$ defined by $Bv=\partial_t (k\ast
v)$, see Section \ref{SecP}.

The paper is organized as follows. In Section 2 we collect some
preliminary results such as the basic inequality (\ref{BI}) and we
explain the time regularization method in more detail. The main
result is proved in Sections 3 and 4. Section 3 is devoted to the
truncated energy estimates and in Section 4 we carry out the
iteration process. Section 5 gives the proof of Theorem
\ref{plaplacebound}. In Section 6 we establish the maximum principle
for homogeneous structures, while Section 7 is concerned with the
case of natural growth conditions.
\section{Preliminaries} \label{SecP}
We first discuss an important method of regularizing kernels of type
${\cal PC}$. Let $k, l \in L_{1,\,loc}(\iR_+)$ be as in assumption
(K1). For $1\le p<\infty$, $T>0$, and a real Banach space $X$ we
consider the operator $B$ defined by
\[ B u=\,\frac{d}{dt}\,(k\ast u),\;\;D(B)=\{u\in L_p([0,T];X):\,k\ast u\in \mbox{}_0 H^1_p([0,T];X)\},
\]
where the zero means vanishing at $t=0$. It is known that this
operator is $m$-accretive in $L_p([0,T];X)$, cf.\ \cite{Phil1},
\cite{CP}, \cite{Grip1}. Its Yosida approximations $B_{n}$, defined
by $B_{n}=nB(n+B)^{-1},\,n\in \iN$, enjoy the property that for any
$u\in D(B)$, one has $B_{n}u\rightarrow Bu$ in $L_p([0,T];X)$ as
$n\to \infty$. It has been shown in \cite{Za} that
\[
B_n u=\,\frac{d}{dt}\,(k_n\ast u),\quad u\in L_p([0,T];X),\;n\in
\iN,
\]
where the kernel $k_n$ has the representation
\begin{equation} \label{knprop}
k_n=k\ast h_n,\quad n\in \iN.
\end{equation}
Here $h_n\in L_{1,\,loc}(\iR_+)$ denotes the resolvent kernel
associated with $nl$, that is
\[
h_n(t)+n(h_n\ast l)(t)=nl(t),\quad t>0,\;n\in\iN.
\]
It is further known that the kernels $k_n$, $n\in\iN$, are also
nonnegative and nonincreasing, and that in addition they belong to
$H^1_1([0,T])$, see e.g.\ \cite{VZ}, \cite{Za}, \cite{ZWH}.

We also remark that (K1) implies that $l$ is completely positive,
see e.g. Theorem 2.2 in \cite{CN}. Consequently, $l$ and $h_n$ are
nonnegative for all $n\in\iN$.

Note further that for any function $f\in L_p([0,T];X)$, $1\le
r<\infty$, there holds $h_n\ast f\to f$ in $L_p([0,T];X)$ as $n\to
\infty$. In fact, setting $u=l\ast f$, we have $u\in D(B)$, and
\[
B_n u=\,\frac{d}{dt}\,(k_n\ast u)=\,\frac{d}{dt}\,(k\ast l\ast
h_n\ast f)=h_n\ast f\,\to\,Bu=f\quad\mbox{in}\;L_p([0,T];X)
\]
as $n\to \infty$. In particular, $k_n\to k$ in $L_1([0,T])$ as $n\to
\infty$.

We next recall a fundamental identity for integro-differential
operators of the form $\frac{d}{dt}(k\ast u)$. Suppose $k\in
H^1_1([0,T])$ and $H\in C^1(\iR)$. Then for any sufficiently smooth
function $u$ on $(0,T)$ one has for a.a. $t\in (0,T)$,
\begin{align} \label{fundidentity}
H'(u(t))&\frac{d}{dt}\,(k \ast u)(t) =\;\frac{d}{dt}\,(k\ast H(u))(t)+
\Big(-H(u(t))+H'(u(t))u(t)\Big)k(t) \nonumber\\
 & +\int_0^t
\Big(H(u(t-s))-H(u(t))-H'(u(t))[u(t-s)-u(t)]\Big)[-\dot{k}(s)]\,ds.
\end{align}
This follows from a straightforward computation, see also \cite{Za},
\cite{Za2}. An integrated version of (\ref{fundidentity}) can be
found in \cite[Lemma 18.4.1]{GLS}.

Equation (\ref{fundidentity}) is highly important for deriving a
priori estimates for problems of the form (\ref{MProb}). In this
paper (see also \cite{Za}) we will apply it to the functions
$H_+(y)=\frac{1}{2}(y_+)^2$ and $H_-(y)=\frac{1}{2}(y_-)^2$ defined
for $y\in \iR$. Here $y_-:=\min\{y,0\}$. Evidently, $H_\pm\in
C^1(\iR)$ with derivative $H_\pm'(y)=y_\pm$, $y\in \iR$. If the
kernel $k$ belongs to $H^1_1([0,T])$ and is nonnegative and
nonincreasing, then it follows from (\ref{fundidentity}) and the
convexity of $H_\pm$ that for any function $u\in L_2([0,T])$,
\begin{equation} \label{BI}
u(t)_\pm\frac{d}{dt}\,(k \ast u)(t)\ge
\frac{1}{2}\,\frac{d}{dt}\,\Big(k \ast (u_\pm)^2\Big)(t),\quad
\mbox{a.a.}\;t\in (0,T).
\end{equation}

The following two lemmas concerning the geometric convergence of
sequences of numbers will be needed for the De Giorgi iteration
arguments below. The first one is contained, e.g., in \cite[Chapter
II, Lemma 5.6]{LSU}, see also \cite[Chapter I, Lemma 4.1]{DB}. Its
proof is by induction.
\begin{lemma} \label{ALG1}
Let $\{Y_n\}$, $n=0,1,2,\ldots$, be a sequence of positive numbers,
satisfying the recursion inequality
\[
Y_{n+1}\le C b^n Y_n^{1+\alpha},\quad n=0,1,2,\ldots,
\]
where $C,\,b>1$ and $\alpha>0$ are given numbers. If
\[
Y_0\le C^{-1/\alpha}b^{-1/\alpha^2},
\]
then
\[
Y_n\le C^{-1/\alpha}b^{-1/\alpha^2}b^{-n/\alpha},\quad n\in \iN,
\]
in particular $Y_n\to 0$ as $n\to \infty$.
\end{lemma}
\begin{lemma} \label{ALG2}
Let $\{Y_n\}$, $n=0,1,2,\ldots$, be a sequence of positive numbers,
satisfying the recursion inequality
\[
Y_{n+1}\le C b^n \big(Y_n^{1+\alpha}+Y_n^{1+\delta}\big),\quad
n=0,1,2,\ldots,
\]
where $C,\,b>1$ and $\delta\ge \alpha>0$ are given numbers. If
\[
Y_0\le (2C)^{-1/\alpha}b^{-1/\alpha^2},
\]
then
\[
Y_n\le (2C)^{-1/\alpha}b^{-1/\alpha^2}b^{-n/\alpha},\quad n\in \iN,
\]
and thus $Y_n\to 0$ as $n\to \infty$.
\end{lemma}
{\em Proof.} The assertion follows directly from the proof of the
previous lemma and the trivial estimate
\[
Y_{n+1}\le 2C b^n Y_n^{1+\alpha},
\]
which holds whenever $Y_n\le 1$, due to the assumption $\delta\ge
\alpha$.\hfill $\square$

$\mbox{}$

We conclude this preliminary part with an interpolation result which
will be frequently used in this paper.

Let $T>0$ and $\Omega$ be a bounded domain in $\iR^N$. For $p,q\ge
1$ we define the spaces
\begin{equation} \label{Vdef}
V_{q,p}:=V_{q,p}([0,T]\times \Omega)=L_{2q}([0,T];L_2(\Omega))\cap
L_p([0,T];H^1_p(\Omega)),
\end{equation}
and
\begin{equation} \label{V0def}
V^0_{q,p}:=V^0_{q,p}([0,T]\times
\Omega)=L_{2q}([0,T];L_2(\Omega))\cap L_p([0,T];\oH^1_p(\Omega)),
\end{equation}
both equipped with the norm
\[
|u|_{V_{q,p}([0,T]\times \Omega)}:=|u|_{L_{2q}([0,T];L_2(\Omega))}
+|Du|_{L_p([0,T];H^1_p(\Omega))}.
\]
We will assume that $\partial \Omega$ satisfies the property of
positive density, i.e. there exist $\delta\in (0,1)$ and $\rho_0>0$
such that for any $x_0\in \Gamma$, any ball $B(x_0,\rho)$ with
$\rho\le \rho_0$ we have that $|\Omega\cap B(x_0,\rho)|\le
\delta|B(x_0,\rho)|$, cf.\ e.g.\ \cite[Section I.1]{DB}.
\begin{prop} \label{parabolicembed}
There exists a constant $\tilde{C}=\tilde{C}(N,p,q)$ such that for
every $u\in V^0_{q,p}([0,T]\times \Omega)$ there holds
\begin{align}
\int_0^T & \int_\Omega |u(t,x)|^r\,dx\,dt \nonumber\\
& \le \tilde{C}^r\Big(\int_0^T \int_\Omega
|Du(t,x)|^p\,dx\,dt\Big)^\frac{\beta r}{p}\,\Big(\int_0^T
\big(\int_\Omega
|u(t,x)|^2\,dx\big)^q\,dt\Big)^\frac{(1-\beta)r}{2q},
\label{multiembed}
\end{align}
where
\begin{equation} \label{rdef}
r=\frac{1-\frac{1}{q}+\frac{2}{N}}{\frac{1}{p}\big(1-\frac{1}{q}\big)+\frac{1}{Nq}}\quad
\mbox{and}\quad\beta=\frac{1-\frac{1}{q}}{1-\frac{1}{q}+\frac{2}{N}}\in
[0,1] .
\end{equation}
\end{prop}
{\em Proof.} We proceed similar as in \cite[Chapter I]{DB}). We
first consider the case where $p>2N/(N+2)$. By the
Gagliardo-Nirenberg inequality (see e.g.\ \cite[Theorem 2.1]{DB}),
we have for a.a. $t\in (0,T)$
\begin{equation} \label{embed1}
|u(t,\cdot)|_{L_r(\Omega)}\le
C_1(N,p)|Du(t,\cdot)|_{L_p(\Omega)}^\beta
|u(t,\cdot)|_{L_2(\Omega)}^{1-\beta},
\end{equation}
where $\beta$ and $r$ are given by (\ref{rdef}). In fact, a short
computation shows that (\ref{rdef}) implies that
\[
\beta=\frac{\frac{1}{2}-\frac{1}{r}}{\frac{1}{N}-\frac{1}{p}+\frac{1}{2}},
\]
which corresponds to condition (2.2) in \cite[Theorem 2.1]{DB}.
Taking the $r$th power in (\ref{embed1}), integrating over $(0,T)$,
and using H\"older's inequality yields
\begin{align*}
|u|_{L_r(\Omega_T)}^r\le C_1^r |Du|_{L_p(\Omega_T)}^{\beta
r}\,|u|_{L_{\hat{r}}([0,T];L_2(\Omega))}^{(1-\beta)r},\quad
\hat{r}=\frac{(1-\beta)rp}{p-\beta r}.
\end{align*}
One verifies that $\hat{r}=2q$, and so (\ref{multiembed}) is valid.

If $p\le 2N/(N+2)$, then we have in particular $p<N$, and thus by
\cite[Corollary 2.1]{DB}
\[
|u(t,\cdot)|_{L_\frac{Np}{N-p}(\Omega)}\le
C_2(N,p)|Du(t,\cdot)|_{L_p(\Omega)},\quad \mbox{a.a.}\;\,t\in (0,T).
\]
Using this and the fact that
\[
r(1-\beta)\,\frac{Np}{Np-(N-p)r\beta}\,=2
\]
it follows by means of H\"older's inequality that
\begin{align*}
\int_0^T \int_\Omega |u|^r\,dx\,dt & =\int_0^T\int_\Omega |u|^{\beta
r}|u|^{(1-\beta)r}\,dx\,dt\\
& \le \int_0^T \big(\int_\Omega
|u|^\frac{Np}{N-p}\,dx\big)^\frac{\beta r(N-p)}{Np}\big(\int_\Omega
|u|^2\,dx\big)^\frac{(1-\beta)r}{2}\,dt\\
& \le C_2^{\beta r}\int_0^T \big(\int_\Omega
|Du|^p\,dx\big)^\frac{\beta r}{p}\big(\int_\Omega
|u|^2\,dx\big)^\frac{(1-\beta)r}{2}\,dt.
\end{align*}
As in the first case we may now apply H\"older's inequality once
more thereby proving (\ref{multiembed}). \hfill $\square$

$\mbox{}$

We remark that $r\ge 2$ if and only if $p\ge \frac{2N}{N+2}$.
\section{Energy estimates} \label{SecEnergy}
The following lemma will be the starting point for all of the a
priori estimates derived in this paper. It provides an equivalent
weak formulation of (\ref{MProb}) where the kernel $k$ is replaced
with the more regular kernel $k_n$ ($n\in \iN$) defined in
(\ref{knprop}). In what follows the kernels $h_n$, $n\in \iN$, are
as in Section \ref{SecP}.
\begin{lemma} \label{LemmaReg}
Let $p>1$, $T>0$ and $\Omega\subset\iR^N$ be a bounded domain. Let
the assumptions (K1),(K2),(Q1)-(Q5) be satisfied and assume that
$u_0\in L_2(\Omega)$. Then $u\in \tilde{V}_{q,p}$ is a weak solution
(subsolution, supersolution) of (\ref{MProb}) if and only if for
every nonnegative function $\psi\in \oH^1_p(\Omega)$ one has
\begin{align}
\int_\Omega \Big(\psi\partial_t [k_n\ast (u-u_0)]+& \big(h_n\ast
a(\cdot,x,u,Du)|D \psi\big)-[h_n\ast
b(\cdot,x,u,Du)]\psi\Big)\,dx\nonumber\\
& =(\le,\,\ge)\,\,0 \quad\mbox{a.a.}\;t\in (0,T),\,n\in \iN.
\label{LemmaRegF}
\end{align}
\end{lemma}
{\em Proof.} The proof is analogous to the proof of Lemma 3.1 in
\cite{Za}. For the reader's convenience we repeat it here.

We may restrict ourselves to the subsolution case as the remaining
cases can be treated analogously.

The 'if' part can be seen as follows. Given an arbitrary nonnegative
$\eta\in \oH^{1,1}_2(\Omega_T)$ satisfying $\eta|_{t=T}=0$, we take
in (\ref{LemmaRegF}) $\psi(x)=\eta(t,x)$ for any fixed $t\in (0,T)$,
integrate from $t=0$ to $t=T$, and integrate by parts w.r.t. the
time variable. Sending then $n\to \infty$ yields (\ref{WQ}); here we
use the approximating properties of the kernels $h_n$ described in
Section \ref{SecP}.

To prove the 'only--if' part, we take the test function
\begin{equation} \label{GB01}
\eta(t,x)=\int_t^T h_n(\sigma-t)\varphi(\sigma,x)\,d\sigma=
\int_0^{T-t}h_n(\sigma)\varphi(\sigma+t,x)\,d\sigma,\quad
t\in(0,T),\,x\in\Omega,
\end{equation}
with arbitrary $n\in \iN$ and nonnegative $\varphi\in
\oH^{1,1}_2(\Omega_T)$ satisfying $\varphi|_{t=T}=0$; $\eta$ is
nonnegative since $\varphi$ and $h_n$ are so (see Section
\ref{SecP}). Then we have
\[
\eta_t(t,x)=\int_t^T
h_n(\sigma-t)\varphi_\sigma(\sigma,x)\,d\sigma,\quad\mbox{a.a.}
\;(t,x)\in \Omega_T.
\]
By Fubini's theorem, we have
\[
\int_{0}^T \Big(\int_t^T h_n(\sigma-t)
\psi_1(\sigma)\,d\sigma\Big)\psi_2(t)\,dt= \int_{0}^T
\psi_1(t)\Big(\int_{0}^t
h_n(t-\sigma)\psi_2(\sigma)\,d\sigma\Big)\,dt,
\]
for all $\psi_1,\,\psi_2\in L_2([0,T])$. So it follows from
(\ref{WQ}) and $k_n=h_n\ast k$ (c.p. (\ref{knprop})) that
\[
\int_0^T\int_\Omega \Big(-\varphi_t[k_n\ast (u-u_0)]+\big(h_n\ast
a(\cdot,x,u,Du)|D \varphi\big)-[h_n\ast b(\cdot,x,u,Du)]\varphi
\Big)\,dx\,dt\le 0,
\]
for all $n\in\iN$. Observe that $k_n\ast (u-u_0)\in\mbox{}_0
H^1_2([0,T];L_2(\Omega))$. Therefore, integrating by parts and using
$\varphi|_{t=T}=0$ yields
\begin{equation} \label{RGWF}
\int_0^T\int_\Omega \Big(\varphi \partial_t[k_n\ast (u-u_0)]
+\big(h_n\ast a(\cdot,x,u,Du)|D \varphi\big)-[h_n\ast
b(\cdot,x,u,Du)]\varphi \Big)\,dx\,dt\le 0,
\end{equation}
for all $n\in \iN$ and $\varphi\in \oH^{1,1}_2(\Omega_T)$ with
$\varphi|_{t=T}=0$. By means of a simple approximation argument, we
then see that (\ref{RGWF}) is valid for any $\varphi$ of the form
$\varphi(t,x)=\chi_{(t_1,t_2)}(t)\psi(x)$, where $\chi_{(t_1,t_2)}$
denotes the characteristic function of the time-interval
$(t_1,t_2)$, $0<t_1<t_2<T$, and $\psi\in \oH^1_2(\Omega)$ is
nonnegative. Relation (\ref{LemmaRegF}) follows now from the
Lebesgue differentiation theorem. \hfill $\square$

$\mbox{}$

Our proof of the sup-bounds for subsolutions stated above relies on
the subsequent truncated energy estimates.
\begin{prop} \label{basicenergy}
Let $p>1$, $T>0$ and $\Omega\subset\iR^N$ be a bounded domain. Let
the assumptions (K1),(K2),(Q1)-(Q5) be satisfied. Suppose that
$u_0\in L_2(\Omega)$ is essentially bounded above in $\Omega$. Then
for any weak subsolution $u\in \tilde{V}_{q,p}$ of (\ref{MProb})
with $\esup_{\Gamma_T} u<\infty$ and any $\kappa$ satisfying the
condition
\begin{equation} \label{GBkcond0}
\kappa\ge
\tilde{\kappa}:=\max\{0,\esup_{\Omega}u_0,\esup_{\Gamma_T} u\},
\end{equation}
there holds
\begin{align*}
|(u-\kappa)_+|^2_{L_{2q}([0,T];L_2(\Omega))}+&\,
|D(u-\kappa)_+|^p_{L_p([0,T]\times \Omega)}\nonumber\\
& \le C\Big(\int_0^T \int_{A_\kappa(t)}u^\gamma
\,dx\,dt+\big(\int_0^T |A_\kappa(t)|\,dt\big)^\frac{s-1}{s}\Big)
\end{align*}
where
\[
A_\kappa(t)=\{x\in \Omega:u(t,x)>\kappa\},\quad t\in (0,T),
\]
and the constant $C=C
(N,p,q,C_0,c_0,C_2,c_2,\gamma,s,|l|_{L_q([0,T])},
|\varphi_0+\varphi_2|_{L_s(\Omega_T)},T,|\Omega|)$.
\end{prop}
{\em Proof.} Let $u\in \tilde{V}_{q,p}$ be a weak subsolution of
(\ref{MProb}) in $\Omega_T$. Then (\ref{LemmaRegF}) holds with the
'$\le$' sign for any nonnegative function $\psi\in \oH^1_p(\Omega)$.
For $t\in (0,T)$ we choose in (\ref{LemmaRegF}) the test function
$\psi=u_{\kappa}^{+}:=(u_\kappa)_+$, where we set
$u_\kappa:=u-\kappa$, and $\kappa\in \iR$ satisfies
(\ref{GBkcond0}). The resulting inequality can be written as
\begin{align}
\int_{\Omega} \Big(u_\kappa^+& \partial_t(k_n\ast
u_\kappa)+\big(h_n\ast a(\cdot,x,u,Du)|D u_\kappa^+\big)\Big)
\,dx\nonumber\\
&\le\int_\Omega \Big([h_n\ast b(\cdot,x,u,Du)]u_\kappa^+
+u_\kappa^+(u_0-\kappa)k_n\Big)\,dx,\quad\mbox{a.a.}\;t\in (0,T).
\label{GB1}
\end{align}
By positivity of $k_n$ and (\ref{GBkcond0}),
\[
\int_{\Omega} u_\kappa^+(u_0-\kappa)k_n\,dx\le
0,\quad\mbox{a.a.}\;t\in (0,T),
\]
Thanks to (\ref{BI}) we further have
\begin{equation} \label{GBFE}
u_\kappa^+\partial_t(k_n\ast u_\kappa)\ge
\,\frac{1}{2}\,\partial_t\Big(k_n\ast (u_\kappa^+)^2\Big),
\quad\mbox{a.a.}\;(t,x)\in\Omega_T.
\end{equation}
Using these relations we infer from (\ref{GB1}) that for a.a.\ $t\in
(0,T)$
\begin{align}
\int_{\Omega} \Big(\frac{1}{2}\,\partial_t [k_n\ast
(u_\kappa^+)^2]+\big(h_n\ast a(\cdot,x,u,Du)|D u_\kappa^+\big)\Big)
\,dx \le \int_\Omega \Big([h_n\ast b(\cdot,x,u,Du)]u_\kappa^+
\,dx.\label{GB2}
\end{align}
We next convolve (\ref{GB2}) with the nonnegative kernel $l$ from
assumption (K1), and observe that in view of
\[
k_n\ast (u_\kappa^+)^2\in \mbox{}_0 H^1_1([0,T];L_1(\Omega))
\]
and $k_n=k\ast h_n$ we have
\[
l\ast\partial_t\Big(k_n\ast
(u_\kappa^+)^2\Big)=\partial_t\Big(l\ast k_n\ast
(u_\kappa^+)^2\Big)=h_n\ast (u_\kappa^+)^2.
\]
Sending then $n\to \infty$, and selecting an appropriate
subsequence, if necessary, we thus obtain
\begin{equation} \label{GB3}
\frac{1}{2}\,\int_{\Omega} (u_\kappa^+)^2\,dx+l\ast \int_{\Omega}
\big(a(\cdot,x,u,Du)|D u_\kappa^+\big)\,dx\le l\ast \int_{\Omega}
b(\cdot,x,u,Du)]u_\kappa^+\,dx
\end{equation}
for a.a.\ $t\in (0,T)$.

By the structure condition (Q1) we have
\begin{align}
\int_{\Omega} \big(a(t,x,u,Du)&|D
u_\kappa^+\big)\,dx=\int_{A_\kappa(t)}\big(a(t,x,u,Du)|D
u\big)\,dx\nonumber\\
& \ge \int_{A_\kappa(t)}
\Big(C_0|Du|^p-c_0|u|^\gamma-\varphi_0\Big)\,dx. \label{GB4}
\end{align}
Employing (Q3) and Young's inequality we may further estimate
\begin{align}
\int_{\Omega} |b(t,x,u,Du) u_\kappa^+|\,dx & \le \int_\Omega
\Big(C_2|Du|^{p\frac{\gamma-1}{\gamma}}u_\kappa^++c_2|u|^{\gamma-1}u_\kappa^++\varphi_2u_\kappa^+\Big)\,dx\nonumber\\
& \le \int_{A_\kappa(t)}
\Big(\frac{C_0}{2}|Du|^p+C_3|u|^\gamma+\varphi_2u_\kappa^+\Big)\,dx,
\label{GB5}
\end{align}
where the constant $C_3>0$ depends only on $C_0,C_2,c_2$, and
$\gamma$. From (\ref{GB3}), (\ref{GB4}), and (\ref{GB5}) we infer
that for a.a.\ $t\in (0,T)$
\begin{equation} \label{GB6}
\int_{\Omega} (u_\kappa^+)^2\,dx+C_0\,l\ast
\int_{\Omega}|Du_\kappa^+|^p\,dx\le 2l\ast F,
\end{equation}
where
\[
F(t)=
\int_{A_\kappa(t)}\Big((C_3+c_0)|u|^\gamma+\varphi_0+\varphi_2u_\kappa^+\Big)\,dx.
\]

Dropping the second term in (\ref{GB6}), which is nonnegative, and
applying Young's inequality for convolutions yields
\begin{align}
|u_\kappa^+|^2_{L_{2q}([0,T];L_2(\Omega))} &=
|(u_\kappa^+)^2|_{L_{q}([0,T];L_1(\Omega))}\nonumber\\
& \le 2|l|_{L_{q}([0,T])} |F|_{L_1([0,T])} \label{GB7}.
\end{align}

On the other hand, we may also drop the first term in (\ref{GB6}),
convolve the resulting inequality with $k$, and use that $k\ast
l=1$, thereby obtaining
\begin{equation} \label{GB8}
C_0 |D u_\kappa^+|^p_{L_p([0,T];L_p(\Omega))}\le 2|F|_{L_1([0,T])}.
\end{equation}

By H\"older's inequality and assumption (Q5) we have
\begin{equation} \label{GB9}
\int_0^T \int_{A_\kappa(t)}\varphi_0\,dx\,dt\le
|\varphi_0|_{L_s(\Omega_T)}\big(\int_0^T
|A_\kappa(t)|\,dt\big)^\frac{1}{s'}.
\end{equation}

Next, set
\begin{equation} \label{etadef}
\eta:=\Big(\frac{\beta}{p}+\frac{(1-\beta)}{2}\Big)^{-1}=
\frac{1-\frac{1}{q}+\frac{2}{N}}{\frac{1}{p}\big(1-\frac{1}{q}\big)+\frac{1}{N}}\,>1.
\end{equation}
Then the term involving $\varphi_2$ can be estimated as follows,
where we use (Q5), H\"older's and Young's inequality, as well as
Proposition \ref{parabolicembed}.
\begin{align}
 & \int_0^T \int_{A_\kappa(t)}\varphi_2
u_\kappa^+\,dx\,dt\le|\varphi_2|_{L_s(\Omega_T)}|u_\kappa^+|_{L_r(\Omega_T)}
\big(\int_0^T |A_\kappa(t)|\,dt\big)^\frac{r-s'}{rs'}\nonumber\\
& \le \tilde{C} |\varphi_2|_{L_s(\Omega_T)}|D
u_\kappa^+|^\beta_{L_p(\Omega_T)}|u_\kappa^+|^{1-\beta}_{L_{2q}([0,T];L_2(\Omega))}\big(\int_0^T
|A_\kappa(t)|\,dt\big)^\frac{r-s'}{rs'}\nonumber\\
& \le \varepsilon^\eta|D
u_\kappa^+|^{\eta\beta}_{L_p(\Omega_T)}|u_\kappa^+|^{\eta(1-\beta)}_{L_{2q}([0,T];L_2(\Omega))}
+\varepsilon^{-\eta'}\tilde{C}^{\eta'}|\varphi_2|_{L_s(\Omega_T)}^{\eta'}\big(\int_0^T
|A_\kappa(t)|\,dt\big)^\frac{\eta'(r-s')}{rs'}\nonumber\\
& \le \varepsilon^\eta
\big(|u_\kappa^+|^2_{L_{2q}([0,T];L_2(\Omega))}+|D
u_\kappa^+|^p_{L_p(\Omega_T)}\big)+\varepsilon^{-\eta'}\tilde{C}^{\eta'}|\varphi_2|_{L_s(\Omega_T)}^{\eta'}\big(\int_0^T
|A_\kappa(t)|\,dt\big)^\frac{\eta'(r-s')}{rs'}, \label{GB9a}
\end{align}
for every $\varepsilon>0$.

Combining (\ref{GB7})-(\ref{GB9a}) and choosing $\varepsilon$ such
that
\[
2\varepsilon^\eta\Big(|l|_{L_{q}([0,T])}+\frac{1}{C_0}\Big)=\,\frac{1}{2}
\]
yields
\begin{align}
|u_\kappa^+&|^2_{L_{2q}([0,T];L_2(\Omega))}+|D
u_\kappa^+|^p_{L_p([0,T];L_p(\Omega))}\nonumber\\
\le &
\;4\Big(|l|_{L_{q}([0,T])}+\frac{1}{C_0}\Big)\Big\{(C_3+c_0)\int_0^T\int_{A_\kappa(t)}|u|^\gamma\,dx\,dt\nonumber\\
& \quad +|\varphi_0|_{L_s(\Omega_T)}\big(\int_0^T
|A_\kappa(t)|\,dt\big)^\frac{1}{s'}+\varepsilon^{-\eta'}\tilde{C}^{\eta'}|\varphi_2|_{L_s(\Omega_T)}^{\eta'}\big(\int_0^T
|A_\kappa(t)|\,dt\big)^\frac{\eta'(r-s')}{rs'}\Big\}. \label{GB11}
\end{align}

Note that
\begin{equation} \label{GB12}
\frac{1}{s'}<\frac{\eta'(r-s')}{rs'}\;\Leftrightarrow\;\frac{s'}{r}<\frac{1}{\eta}\;\Leftrightarrow\;
s>\frac{r}{r-\eta}=\frac{\frac{1}{p}\big(1-\frac{1}{q}\big)+\frac{1}{N}}{\frac{1}{N}\big(1-\frac{1}{q}\big)}.
\end{equation}
The last condition is exactly the one in (Q5). Hence we may estimate
\[
\big(\int_0^T |A_\kappa(t)|\,dt\big)^\frac{\eta'(r-s')}{rs'}\le
\big(T|\Omega|\big)^{\frac{\eta'(r-s')}{rs'}-\frac{1}{s'}}\big(\int_0^T
|A_\kappa(t)|\,dt\big)^\frac{1}{s'},
\]
and thus (\ref{GB11}) implies the desired energy estimate. \hfill
$\square$

$\mbox{}$

The corresponding result for supersolutions reads as follows.
\begin{prop} \label{basicenergy2}
Let $p>1$, $T>0$ and $\Omega\subset\iR^N$ be a bounded domain. Let
the assumptions (K1),(K2),(Q1)-(Q5) be satisfied. Suppose that
$u_0\in L_2(\Omega)$ is essentially bounded below in $\Omega$. Then
for any weak supersolution $u\in \tilde{V}_{q,p}$ of (\ref{MProb})
with $\einf_{\Gamma_T} u>-\infty$ and any $\kappa$ satisfying the
condition
\begin{equation} \label{GBkcond}
\kappa\ge
\hat{\kappa}:=-\min\{0,\einf_{\Omega}u_0,\einf_{\Gamma_T} u\},
\end{equation}
we have
\begin{align}
|(u+\kappa)_-|^2_{L_{2q}([0,T];L_2(\Omega))}+&\,
|D(u+\kappa)_-|^p_{L_p([0,T]\times \Omega)}\nonumber\\
& \le C\Big(\int_0^T \int_{\tilde{A}_\kappa(t)}(-u)^\gamma
\,dx\,dt+\big(\int_0^T
|\tilde{A}_\kappa(t)|\,dt\big)^\frac{s-1}{s}\Big)
\end{align}
where
\[
\tilde{A}_\kappa(t)=\{x\in \Omega:-u(t,x)>\kappa\},\quad t\in
(0,T),
\]
and the constant $C$ is like in Proposition \ref{basicenergy}.
\end{prop}
{\em Proof.} The proof is analogous to the previous one.
(\ref{LemmaRegF}) now holds with the '$\ge$' sign and we take
$\psi=-(u+\kappa)_-\ge 0$. Replacing $u$ by $-u$ and
${A}_\kappa(t)$ by $\tilde{A}_\kappa(t)$ the same line of
arguments as above yields the asserted estimate. \hfill $\square$
\section{Iterative inequalities} \label{SecIterations}
Let $u\in \tilde{V}_{q,p}$ be a weak subsolution of (\ref{MProb}) in
$\Omega_T$. Set
\[
\kappa_n=\kappa\Big(2-\frac{1}{2^n}\Big),\quad n=0,1,2,\ldots,
\]
where $\kappa\ge \max\{\tilde{\kappa},1\}$ will be chosen later. We
further put
\[
Y_n=\int_0^T \int_{A_{\kappa_n(t)}}
(u-\kappa_n)_+^\gamma\,dx\,dt,\quad n=0,1,2,\ldots
\]
By Proposition \ref{basicenergy} we have for all $n=0,1,2,\ldots$
\begin{align}
|(u-\kappa_{n+1})_+|^2_{L_{2q}([0,T];L_2(\Omega))}+&\,
|D(u-\kappa_{n+1})_+|^p_{L_p([0,T]\times \Omega)}\nonumber\\
& \le C\Big(\int_0^T \int_{A_{\kappa_{n+1}(t)}}u^\gamma
\,dx\,dt+\big(\int_0^T
|A_{\kappa_{n+1}(t)}|\,dt\big)^\frac{1}{s'}\Big). \label{itera1}
\end{align}
To estimate the right-hand side, note first that
\begin{align}
\int_0^T |A_{\kappa_{n+1}(t)}| & \,dt \le \int_0^T
\int_{A_{\kappa_{n+1}(t)}}\Big(\frac{u-\kappa_n}{\kappa_{n+1}-\kappa_n}\Big)^\gamma
\,dx\,dt \nonumber\\
& \le \frac{1}{(\kappa_{n+1}-\kappa_n)^\gamma}\int_0^T
\int_{A_{\kappa_{n}(t)}}(u-\kappa_n)^\gamma\,dx\,dt =
\,\frac{2^{\gamma(n+1)}}{\kappa^\gamma}\,Y_n. \label{itera2}
\end{align}
Further,
\begin{align*}
Y_n \ge \int_0^T
\int_{A_{\kappa_{n+1}(t)}}&\,(u-\kappa_n)_+^\gamma\,dx\,dt\ge
\int_0^T \int_{A_{\kappa_{n+1}(t)}} u^\gamma
\Big(1-\frac{\kappa_n}{\kappa_{n+1}}\Big)^\gamma\,dx\,dt\\
& \ge \,\frac{1}{2^{\gamma(n+2)}}\int_{A_{\kappa_{n+1}(t)}}
u^\gamma\,dx\,dt.
\end{align*}
Hence (\ref{itera1}) implies that
\begin{align}
|(u-\kappa_{n+1})_+|^2_{L_{2q}([0,T];L_2(\Omega))}+&\,
|D(u-\kappa_{n+1})_+|^p_{L_p([0,T]\times \Omega)}\nonumber\\
& \le
C\,2^{\gamma(n+2)}Y_n+C\,\Big(\frac{2^{\gamma(n+1)}}{\kappa^\gamma}
\Big)^\frac{1}{s'}Y_n^\frac{1}{s'}. \label{itera3}
\end{align}

On the other hand, we have by H\"older's inequality and Proposition
\ref{parabolicembed}
\begin{align*}
Y_{n+1} & =\int_0^T \int_{A_{\kappa_{n+1}(t)}}
(u-\kappa_{n+1})_+^\gamma\,dx\,dt\\
& \le \Big(\int_0^T \int_\Omega
(u-\kappa_{n+1})_+^r\,dx\,dt\Big)^\frac{\gamma}{r} \Big(\int_0^T
|A_{\kappa_{n+1}(t)}|\,dt\Big)^{1-\frac{\gamma}{r}}\\
& \le \tilde{C}^\gamma \Big(|D(u-\kappa_{n+1})_+|_{L_p([0,T]\times
\Omega)}\Big)^{\beta\gamma}
\Big(|(u-\kappa_{n+1})_+|_{L_{2q}([0,T];L_2(\Omega))}\Big)^{(1-\beta)\gamma}\\
& \quad\quad \times \Big(\int_0^T
|A_{\kappa_{n+1}(t)}|\,dt\Big)^{1-\frac{\gamma}{r}},
\end{align*}
where $\tilde{C}=\tilde{C}(N,p,q)$ and $r$ and $\beta$ are given by
(\ref{rdef}). Recall that
\[
\frac{1}{\eta}=\frac{\beta}{p}+\frac{(1-\beta)}{2}=\frac{\frac{1}{p}\big(1-\frac{1}{q}\big)+\frac{1}{N}}{1-\frac{1}{q}+\frac{2}{N}}\,<1.
\]
Using (\ref{itera2}) and (\ref{itera3}) it then follows that
\begin{align*}
Y_{n+1} & \le \tilde{C}^\gamma
\Big[C\,2^{\gamma(n+2)}Y_n+C\,\Big(\frac{2^{\gamma(n+1)}}{\kappa^\gamma}
\Big)^\frac{1}{s'}Y_n^\frac{1}{s'}\Big]^{\frac{\gamma}{\eta}}
\Big[\frac{2^{\gamma(n+1)}}{\kappa^\gamma}\,Y_n\Big]^{1-\frac{\gamma}{r}}\\
& \le \tilde{C}^\gamma
(2C)^{\frac{\gamma}{\eta}}\kappa^{-\gamma(1-\frac{\gamma}{r})}2^{\gamma(n+2)(\frac{\gamma}{\eta}+1)}
\Big(Y_n^{\frac{\gamma}{\eta s'}+1-\frac{\gamma}{r}}
+Y_n^{\frac{\gamma}{\eta}+1-\frac{\gamma}{r}}\Big)\\
& \le \Big(C\tilde{C}4^{\gamma+2}\Big)^\gamma
\kappa^{-\gamma(1-\frac{\gamma}{r})}
\Big(2^{\gamma(\gamma+1)}\Big)^n\Big(Y_n^{\frac{\gamma}{\eta
s'}+1-\frac{\gamma}{r}}
+Y_n^{\frac{\gamma}{\eta}+1-\frac{\gamma}{r}}\Big).
\end{align*}
Note that (see also (\ref{GB12}))
\[
\frac{\gamma}{\eta s'}+1-\frac{\gamma}{r}>1 \;\Leftrightarrow\;
\frac{1}{\eta s'}>\frac{1}{r}\;\Leftrightarrow\; s>\frac{r}{r-\eta
}=\frac{\frac{1}{p}\big(1-\frac{1}{q}\big)+\frac{1}{N}}
{\frac{1}{N}\big(1-\frac{1}{q}\big)}.
\]
The last condition is satisfied thanks to (Q5). Hence
\[
\alpha:=\gamma\Big(\frac{1}{\eta s'}-\frac{1}{r}\Big)>0
\]
as well as
\[
\delta:=\gamma\Big(\frac{1}{\eta}-\frac{1}{r}\Big)\ge \alpha>0.
\]

It follows from Lemma \ref{ALG2} that $Y_n\to 0$ as $n\to \infty$,
provided that
\[
Y_0=\int_0^T \int_{\Omega} (u-\kappa)_+^\gamma\,dx\,dt\le
\Big(\frac{2C^\gamma\tilde{C}^\gamma
4^{\gamma(\gamma+2)}}{\kappa^{\gamma(1-\frac{\gamma}{r})}}\Big)^{-\frac{1}{\alpha}}\big(2^{\gamma(\gamma+1)}\big)^{-\frac{1}{\alpha^2}},
\]
which in turn is certainly satisfied if
\begin{equation} \label{itera5}
\int_0^T \int_{\Omega} u_+^\gamma\,dx\,dt\le \hat{C}^{-1}
\kappa^{\frac{\gamma}{\alpha}(1-\frac{\gamma}{r})},
\end{equation}
where
\[
\hat{C}=\Big(2C^\gamma\tilde{C}^\gamma
4^{\gamma(\gamma+2)}\Big)^{\frac{1}{\alpha}}\big(2^{\gamma(\gamma+1)}\big)^{\frac{1}{\alpha^2}}.
\]
The conditions (\ref{itera5}) and $\kappa\ge
\max\{\tilde{\kappa},1\}$ are fulfilled when we set
\[
\kappa:=\tilde{\kappa}+\max\Big\{1,\big(\hat{C}\int_0^T
\int_{\Omega} u_+^\gamma\,dx\,dt\big)^{\frac{\alpha
r}{\gamma(r-\gamma)}}\Big\}.
\]
Since $\kappa_n\to 2\kappa$ as $n\to \infty$ we thus obtain
\[
\esup_{\Omega_T} u\le
2\kappa=2\Big(\tilde{\kappa}+\max\Big\{1,\big(\hat{C}\int_0^T
\int_{\Omega} u_+^\gamma\,dx\,dt\big)^{\frac{\alpha
r}{\gamma(r-\gamma)}}\Big\}\Big).
\]

Finally a short computation shows that
\[
\frac{\alpha
r}{\gamma}=\frac{\frac{1}{N}\big(1-\frac{1}{q}\big)-\frac{1}{s}\big[\frac{1}{p}\big(1-\frac{1}{q}\big)+\frac{1}{N}\big]}
{\frac{1}{p}\big(1-\frac{1}{q}\big)+\frac{1}{Nq}}.
\]
The first part of the theorem is proved.

The second part is proved analogously replacing $u$ by $-u$ and
${A}_\kappa(t)$ by $\tilde{A}_\kappa(t)$ in the previous arguments
and employing Proposition \ref{basicenergy2}. \hfill $\square$
\section{Proof of Theorem \ref{plaplacebound}} \label{pl}
To prove Theorem \ref{plaplacebound} note first that the kernel
$k=g_{1-\alpha}$ satisfies (K1) and (K2) with $l=g_\alpha$ and any
$q\in (1,\frac{1}{1-\alpha})$. In particular, any $q>1$ satisfying
condition (\ref{qcond}) is admissible.

Let now $q>1$ be fixed such that (\ref{qcond}) holds and suppose
that $u\in \tilde{V}_{q,p}$ is a weak solution of (\ref{plaplace})
in $\Omega_T$. Fix $\gamma\in (1,\eta)$, where $\eta$ is given by
(\ref{etadef}). Note that $\eta<r$. Hence we may apply Theorem
\ref{mainresult}, which yields boundedness of $u$ in $\Omega_T$ and
the estimate
\begin{equation} \label{mainbound1}
|u|_{L_\infty(\Omega_T)}\le
2\Big(\max\{|u_0|_{L_\infty(\Omega)},\esup_{\Gamma_T}
|u|\}+\max\Big\{1,C\big(\int_0^T \int_{\Omega}
|u|^\gamma\,dx\,dt\big)^{\frac{\theta}{r-\gamma}}\Big\}\Big),
\end{equation}
where the constant $C$ depends only on the data.

It remains to derive an a priori bound for the integral term
involving $|u|$ in terms of the data and the quantity
$\kappa:=\max\{|u_0|_{L_\infty(\Omega)},\esup_{\Gamma_T} |u|\}$. To
this purpose we write $|u|=u_++(-u)_+$ and estimate
$|u_+|_{L_\gamma(\Omega_T)}$ and $|(-u)_+|_{L_\gamma(\Omega_T)}$
separately by means of the energy estimates from Propositions
\ref{basicenergy} and \ref{basicenergy2}, respectively.

We have
\begin{align*}
\int_0^T \int_\Omega u_+^\gamma\,dx\,dt\le 2^\gamma \Big(\int_0^T
\int_\Omega (u-\kappa)_+^\gamma\,dx\,dt+\kappa^\gamma
T|\Omega|\Big),
\end{align*}
and by Propositions \ref{parabolicembed} and \ref{basicenergy}
\begin{align*}
\int_0^T & \int_\Omega (u-\kappa)_+^\gamma\,dx\,dt \le
|(u-\kappa)_+|_{L_r(\Omega_T)}^\gamma
\big(T|\Omega|\big)^\frac{r-\gamma}{\gamma}\\
& \le
\varepsilon^{\frac{\eta}{\gamma}}|(u-\kappa)_+|_{L_r(\Omega_T)}^{\eta}+\varepsilon^{-\frac{\eta}{\eta-\gamma}}
\big(T|\Omega|\big)^{\frac{(r-\gamma)\eta}{\gamma(\eta-\gamma)}}\\
& \le \varepsilon^{\frac{\eta}{\gamma}} \tilde{C}^\eta
\Big(|(u-\kappa)_+|^2_{L_{2q}([0,T];L_2(\Omega))}+
|D(u-\kappa)_+|^p_{L_p([0,T]\times
\Omega)}\Big)+\varepsilon^{-\frac{\eta}{\eta-\gamma}}
\big(T|\Omega|\big)^{\frac{(r-\gamma)\eta}{\gamma(\eta-\gamma)}}\\
& \le \varepsilon^{\frac{\eta}{\gamma}} \tilde{C}^\eta
C\Big(\int_0^T \int_\Omega
u_+^\gamma\,dx\,dt+\big(T|\Omega|)^\frac{1}{s'}\Big)+\varepsilon^{-\frac{\eta}{\eta-\gamma}}
\big(T|\Omega|\big)^{\frac{(r-\gamma)\eta}{\gamma(\eta-\gamma)}},
\end{align*}
for all $\varepsilon>0$. Choosing $\varepsilon$ sufficiently small
we get a bound
\[
\int_0^T \int_\Omega u_+^\gamma\,dx\,dt\le
C=C(N,p,q,s,|f|_{L_s(\Omega_T)},T,|\Omega|,\kappa).
\]
The bound for $|(-u)_+|_{L_\gamma(\Omega_T)}$ is obtained
analogously. Combining these estimates and (\ref{mainbound1}) proves
the assertion of Theorem \ref{plaplacebound}.  \hfill $\square$

$\mbox{}$

The condition on $f$ is sharp, at least in the cases $p=2$,
$\alpha\in(0,1)$, and $p>2$, $\alpha\in(p'/N,1)$, as we will show in
the following.

Suppose $u$ is a solution of equation (\ref{plaplace}) with smooth
data $u_0$, e.g.\ $u_0=0$. Then the optimal regularity for $u$ is
determined from the conditions
\begin{equation} \label{MRcond}
\partial_t^\alpha (u-u_0)\in L_s(\Omega_T)\quad \mbox{and}\;\;\,\mbox{div}\,\big(|Du|^{p-2}Du\big)\in
L_s(\Omega_T).
\end{equation}
In the linear case $p=2$ these lead to the maximal regularity class
\[
u\in Z:=H^\alpha_s([0,T];L_s(\Omega))\cap
L_{s}([0,T];H^2_{s}(\Omega)),
\]
where $H^\alpha_s(J;X)$ denotes the vector-valued Bessel potential
space of $X$-valued functions on the interval $J$. In fact, the
first space is a consequence of the first condition, by well-known
properties of the fractional derivation operator, see e.g.\
\cite{ZEQ}.

The question is now for which $s>1$ we have $Z\hookrightarrow
L_\infty(\Omega_T)$. This can be determined by means of cross
interpolation and Sobolev embeddings. By the mixed derivative
theorem (see \cite{Sob}) we have for all $\vartheta\in [0,1]$
\[
Z\hookrightarrow
H^{(1-\vartheta)\alpha}_s([0,T];H^{2\vartheta}_s(\Omega)).
\]
Thus $Z\hookrightarrow L_\infty(\Omega_T)$ if there exists
$\vartheta\in(0,1)$ such that
\[
(1-\vartheta)\alpha>\frac{1}{s}\quad\mbox{and}\quad
2\vartheta>\frac{N}{s}.
\]
This is equivalent to
\begin{equation} \label{bed1}
\frac{N}{2s}<\vartheta<1-\frac{1}{s\alpha}.
\end{equation}
There exists $\vartheta\in [0,1]$ satisfying (\ref{bed1}) if and
only if $s>\frac{N}{2}+\frac{1}{\alpha}$, which is exactly the
condition for $s$ in Theorem \ref{plaplacebound} in the case $p=2$.

We now discuss a nonlinear case, namely let us assume that $p>2$ and
$\alpha\in(p'/N,1)$. Here we are led to the maximal regularity class
\[
u\in Z:=H^\alpha_s([0,T];L_s(\Omega))\cap
L_{s(p-1)}([0,T];H^2_{\hat{r}}(\Omega))\quad
\mbox{with}\;\;\frac{p-1}{\hat{r}}=\frac{1}{s}+\frac{p-2}{N}.
\]
To see the second space suppose that $u\in
L_{\bar{q}}([0,T];H^2_{\bar{r}}(\Omega))$ satisfies the second
condition in (\ref{MRcond}). Assuming $\bar{r}<N$ we have $D_i D_j
u\in L_{\bar{q}}([0,T];L_{\bar{r}}(\Omega))$ and
\[
D_i u\in L_{\bar{q}}([0,T];H^1_{\bar{r}}(\Omega))\hookrightarrow
L_{\bar{q}}([0,T];L_{\frac{\bar{r}N}{N-\bar{r}}}(\Omega)),
\]
for all $i,j=1,\ldots,n$, and so the structure of the $p$-Laplacian
leads to the condition
\[
L_{\frac{\bar{q}}{p-2}}([0,T];L_{\frac{\bar{r}N}{(N-\bar{r})(p-2)}}(\Omega))\times
L_{\bar{q}}([0,T];L_{\bar{r}}(\Omega))\subset
L_s([0,T];L_s(\Omega)).
\]
This means we need
\[
\frac{p-2}{\bar{q}}+\frac{1}{\bar{q}}=\frac{1}{s}\quad
\mbox{and}\;\;\,\frac{(N-\bar{r})(p-2)}{\bar{r}N}+\frac{1}{\bar{r}}=\frac{1}{s},
\]
which in turn implies $\bar{q}=s(p-1)$ and $\bar{r}=\hat{r}$.

As before we have to determine those $s>1$ for which we have
$Z\hookrightarrow L_\infty(\Omega_T)$. Note first that the
conditions $\alpha>\frac{1}{s}$ and $2>\frac{N}{\hat{r}}$ are
necessary. Assuming this it follows that
$s>\max\{\frac{1}{\alpha},\frac{N}{p}\}$, in particular we have
$s\ge\hat{r}$ due to $p> 2$. By the mixed derivative theorem we then
have for all $\vartheta\in [0,1]$
\begin{align*}
Z & \hookrightarrow
H^{\alpha-\frac{1}{s}+\frac{1}{s(p-1)}}_{s(p-1)}([0,T];L_s(\Omega))\cap
L_{s(p-1)}([0,T];H_s^{2-\frac{N}{\hat{r}}+\frac{N}{s}}(\Omega))\\
& \hookrightarrow
H^{(1-\vartheta)(\alpha-\frac{1}{s}+\frac{1}{s(p-1)})}_{s(p-1)}([0,T];
H_s^{\vartheta(2-\frac{N}{\hat{r}}+\frac{N}{s})}(\Omega))
\end{align*}
with sharp embeddings. Hence $Z\hookrightarrow L_\infty(\Omega_T)$
if there exists $\vartheta\in (0,1)$ such that
\[
(1-\vartheta)\Big(\alpha-\frac{1}{s}+\frac{1}{s(p-1)}\Big)>\frac{1}{s(p-1)}
\quad
\mbox{and}\;\;\;\vartheta\Big(2-\frac{N}{\hat{r}}+\frac{N}{s}\Big)>\frac{N}{s},
\]
which is equivalent to
\[
0<\omega_1:=\frac{\frac{N}{s}}{2-\frac{N}{\hat{r}}+\frac{N}{s}}
<\vartheta<\frac{\alpha-\frac{1}{s}}{\alpha-\frac{1}{s}+\frac{1}{s(p-1)}}=:\omega_2<1.
\]
Thus it boils down to the condition $\omega_1<\omega_2$. A short
computation shows that this condition is in fact equivalent to
$s>\frac{N}{p}+\frac{1}{\alpha}$.

Recall that we assumed that $\bar{r}<N$, that is
$\hat{r}=\bar{r}<N$. The condition $\hat{r}<N$ in turn is equivalent
to $s<N$. Thus $s$ has to satisfy the condition
\[
\frac{N}{p}+\frac{1}{\alpha}<s<N,
\]
which is possible, by the assumption $\alpha\in(p'/N,1)$. Hence in
this case the condition on $f$ is optimal as well.
\section{Homogenous structures} \label{SecHomStruc}
In this section we consider the special case of homogenous
structures. By this we mean equations of the type
\begin{equation} \label{HProb}
\partial_t \Big(k\ast(u-u_0)\Big)-\mbox{div}\,a(t,x,u,Du)=0,\;\;t\in (0,T),\,x\in
\Omega,
\end{equation}
where
\begin{itemize}
\item[{\bf (HS)}] $\quad\quad\quad\quad(a(t,x,\xi,\eta)|\eta) \ge \,
C_0|\eta|^p$,  $\quad\quad\quad|a(t,x,\xi,\eta)| \,\le \,
C_1|\eta|^{p-1}$,
\end{itemize}
for a.a. $(t,x)\in \Omega_T$, and all $\xi\in\iR$, $\eta\in
\iR^N$. Here $C_0$ and $C_1$ are positive constants.

In this situation the weak maximum principle takes the same form
as in the classical parabolic case. Moreover the assumption (K2)
can be dropped, $q=1$ is here admissible.
\begin{satz} \label{MaxPrinciple}
Let $p>1$, $T>0$, and $\Omega\subset\iR^N$ be a bounded domain.
Assume (K1) and (HS), and let $u_0\in L_2(\Omega)$. Then for any
weak subsolution (supersolution) $u\in \tilde{V}_{1,p}$ of
(\ref{HProb}), we have for a.a. $(t,x)\in \Omega_T$
\[
u(t,x)\le \max\Big\{0,\esup_\Omega
u_0,\esup_{\Gamma_T}u\Big\}\quad\quad \Big( \; u(t,x)\ge
\min\Big\{0,\einf_{\Omega}u_0,\einf_{\Gamma_T}u\Big\}\;\Big),
\]
provided this maximum (minimum) is finite.
\end{satz}
{\em Proof.} Note first that Lemma \ref{LemmaReg} also holds under
the assumptions of Theorem \ref{MaxPrinciple}. It suffices to
consider the subsolution case. We take
\begin{equation} \label{kappawahl}
\kappa=\tilde{\kappa}=\max\Big\{0,\esup_\Omega
u_0,\esup_{\Gamma_T}u\Big\},
\end{equation}
assuming that this quantity is finite, and proceed as in the proof
of Proposition \ref{basicenergy}. This yields
\[
|(u-\tilde{\kappa})_+|^2_{L_{2}(\Omega_T)}+
|D(u-\tilde{\kappa})_+|^p_{L_p(\Omega_T)}\le 0,
\]
which immediately implies the assertion. \hfill $\square$

$\mbox{}$

Theorem \ref{MaxPrinciple} shows in particular that the maximum
principle holds in the usual form for weak solutions of the time
fractional $p$-Laplace equation (\ref{plaplace}) with $f=0$.

We remark that the case $q=1$ can occur. In \cite[Section 3]{Za}
an example is given for a kernel $k$ satisfying (K1) with $l\notin
L_{\hat{q}}([0,T])$ for all $\hat{q}>1$ and $T>0$.
\section{Natural growth conditions} \label{SecHGC}
Finally we consider the case of 'natural' or Hadamard growth
conditions with respect to $|Du|$. For the sake of simplicity we
suppose that
\begin{itemize}
\item[{\bf (Q)}] $\quad\quad(a(t,x,\xi,\eta)|\eta) \ge \,
C_0|\eta|^p$,  $\quad|a(t,x,\xi,\eta)| \,\le \, C_1|\eta|^{p-1}$,
$\quad|b(t,x,\xi,\eta)| \,\le \, C_2|\eta|^p$,
\end{itemize}
for a.a. $(t,x)\in \Omega_T$, and all $\xi\in\iR$, $\eta\in
\iR^N$, where $C_i$, $i=0,1,2$ are positive constants. In the
classical parabolic case it is known that weak solutions of the
corresponding problem under the conditions (Q) are in general not
bounded. However there are results (also in a more general
situation) which provide $L_\infty$ bounds in terms of the data
assuming in addition that the weak solution is bounded, see e.g.\
\cite{DB}. It turns out that corresponding results can be obtained
for (\ref{MProb}). Here we only prove such a result in the case
where (Q) holds. It generalizes Theorem 4.3 in \cite{Za}, where
$p=2$ is required. As in the previous section we may drop
assumption (K2).
\begin{satz} \label{MaxPrinciple4}
Let $p>1$, $T>0$, and $\Omega\subset\iR^N$ be a bounded domain.
Let (K1) and (Q) be satisfied, and suppose that $u_0\in
L_\infty(\Omega)$. Then for any bounded weak solution $u\in
\tilde{V}_{1,p}$ of (\ref{MProb}),
\[
|u|_{L_\infty(\Omega_T)}\le
\max\Big\{|u_0|_{L_\infty(\Omega)},\esup_{\Gamma_T}|u|\Big\}.
\]
\end{satz}
{\em Proof.} We proceed as in the proof of \cite[Theorem
17.1]{DB}, see also \cite[Theorem 4.3]{Za}. Assume that
$K:=\esup_{\Omega_T} u>\tilde{\kappa}$, where $\tilde{\kappa}$ is
as in (\ref{kappawahl}). Let $\varepsilon>0$ be such that
$\kappa:=K-\varepsilon\ge \tilde{\kappa}$. We then choose the test
functions $u_\kappa^+=(u-\kappa)_+$ and estimate similarly as
above, using the conditions (Q). This yields
\begin{align*}
|u_\kappa^+|^2_{L_{2}(\Omega_T)}+ |Du_\kappa^+|^p_{L_p(\Omega_T)}
& \le
C(C_0,C_2)\Big||Du_\kappa^+|^p\,u_\kappa^+\Big|_{L_1(\Omega_T)}\\
& \le \varepsilon
C(C_0,C_2)\Big||Du_\kappa^+|^p\Big|_{L_1(\Omega_T)}.
\end{align*}
Choosing $\varepsilon$ sufficiently small, it follows that
$|u_\kappa^+|^2_{L_{2}(\Omega_T)}\le 0$, that is $u\le \kappa<K$
a.e.\ in $\Omega_T$, a contradiction. Hence $u\le \tilde{\kappa}$
a.e.\ in $\Omega_T$. The lower bound is obtained analogously.
 \hfill $\square$

$\mbox{}$

\noindent {\footnotesize {\bf Vicente Vergara}, Universidad de
Tarapac\'{a}, Instituto de Alta Investigaci\'{o}n, Antofagasta N.
1520, Arica, Chile, Email: vvergaraa@uta.cl

$\mbox{}$

\noindent {\bf Rico Zacher}, Martin-Luther-Universit\"at
Halle-Wittenberg, Institut f\"ur Mathematik, Theodor-Lieser-Strasse
5, 06120 Halle, Germany, Email: rico.zacher@mathematik.uni-halle.de

}
\end{document}